%
%
%
%
%

\font\elevensc=cmcsc10 scaled\magstephalf

\font\eightsc=cmcsc10 scaled800

\font\elevenrm=cmr10 scaled \magstephalf
\font\ninerm=cmr9
\font\eightrm=cmr8
\font\sixrm=cmr6
\font\fiverm=cmr5

\font\eleveni=cmmi10 scaled\magstephalf
\font\ninei=cmmi9
\font\eighti=cmmi8
\font\sixi=cmmi6
\font\fivei=cmmi5
\skewchar\ninei='177 \skewchar\eighti='177 \skewchar\sixi='177
\skewchar\eleveni= '177

\font\elevensy=cmsy10 scaled\magstephalf
\font\ninesy=cmsy9
\font\eightsy=cmsy8
\font\sixsy=cmsy6
\font\fivesy=cmsy5
\skewchar\ninesy='60 \skewchar\eightsy='60 \skewchar\sixsy='60
\skewchar\elevensy='60

\font\eighteenbf=cmbx10 scaled\magstep3

\font\twelvebf=cmbx10 scaled \magstep1
\font\elevenbf=cmbx10 scaled \magstephalf
\font\tenbf=cmbx10
\font\ninebf=cmbx9
\font\eightbf=cmbx8
\font\sixbf=cmbx6
\font\fivebf=cmbx5

\font\elevenit=cmti10 scaled\magstephalf
\font\nineit=cmti9
\font\eightit=cmti8

\font\eighteenmib=cmmib10 scaled \magstep3
\font\twelvemib=cmmib10 scaled \magstep1
\font\elevenmib=cmmib10 scaled\magstephalf
\font\tenmib=cmmib10
\font\eightmib=cmmib10 scaled 800
\font\sixmib=cmmib10 scaled 600

\font\eighteensyb=cmbsy10 scaled \magstep3
\font\twelvesyb=cmbsy10 scaled \magstep1
\font\elevensyb=cmbsy10 scaled \magstephalf
\font\tensyb=cmbsy10
\font\eightsyb=cmbsy10 scaled 800
\font\sixsyb=cmbsy10 scaled 600

\font\elevenex=cmex10 scaled \magstephalf
\font\tenex=cmex10
\font\eighteenex=cmex10 scaled \magstep3


\def\elevenpoint{\def\rm{\fam0\elevenrm}%
  \textfont0=\elevenrm \scriptfont0=\eightrm \scriptscriptfont0=\sixrm
  \textfont1=\eleveni \scriptfont1=\eighti \scriptscriptfont1=\sixi
  \textfont2=\elevensy \scriptfont2=\eightsy \scriptscriptfont2=\sixsy
  \textfont3=\elevenex \scriptfont3=\elevenex \scriptscriptfont3=\elevenex
  \def\bf{\fam\bffam\elevenbf}%
  \def\it{\fam\itfam\elevenit}%
  \textfont\bffam=\elevenbf \scriptfont\bffam=\eightbf
   \scriptscriptfont\bffam=\sixbf
\normalbaselineskip=13.9pt
  \setbox\strutbox=\hbox{\vrule height9.5pt depth4.4pt width0pt\relax}%
  \normalbaselines\rm}

\elevenpoint 

\def\ninepoint{\def\rm{\fam0\ninerm}%
  \textfont0=\ninerm \scriptfont0=\sixrm \scriptscriptfont0=\fiverm
  \textfont1=\ninei \scriptfont1=\sixi \scriptscriptfont1=\fivei
  \textfont2=\ninesy \scriptfont2=\sixsy \scriptscriptfont2=\fivesy
  \textfont3=\tenex \scriptfont3=\tenex \scriptscriptfont3=\tenex
  \def\it{\fam\itfam\nineit}%
  \textfont\itfam=\nineit
  \def\bf{\fam\bffam\ninebf}%
  \textfont\bffam=\ninebf \scriptfont\bffam=\sixbf
   \scriptscriptfont\bffam=\fivebf
\normalbaselineskip=11pt
  \setbox\strutbox=\hbox{\vrule height8pt depth3pt width0pt\relax}%
  \normalbaselines\rm}

\def\eightpoint{\def\rm{\fam0\eightrm}%
  \textfont0=\eightrm \scriptfont0=\sixrm \scriptscriptfont0=\fiverm
  \textfont1=\eighti \scriptfont1=\sixi \scriptscriptfont1=\fivei
  \textfont2=\eightsy \scriptfont2=\sixsy \scriptscriptfont2=\fivesy
  \textfont3=\tenex \scriptfont3=\tenex \scriptscriptfont3=\tenex
  \def\it{\fam\itfam\eightit}%
  \textfont\itfam=\eightit
  \def\bf{\fam\bffam\eightbf}%
  \textfont\bffam=\eightbf \scriptfont\bffam=\sixb
   \scriptscriptfont\bffam=\fivebf
\normalbaselineskip=12pt
  \setbox\strutbox=\hbox{\vrule height8.5pt depth3.5pt width0pt\relax}%
  \normalbaselines\rm}


\def\eighteenbold{\def\rm{\fam0\eighteenrm}%
  \textfont0=\eighteenbf \scriptfont0=\twelvebf \scriptscriptfont0=\tenbf
  \textfont1=\eighteenmib \scriptfont1=\twelvemib\scriptscriptfont1=\tenmib
  \textfont2=\eighteensyb \scriptfont2=\twelvesyb\scriptscriptfont2=\tensyb
  \textfont3=\eighteenex \scriptfont3=\tenex \scriptscriptfont3=\tenex
  \def\bf{\fam\bffam\eighteenbf}%
  \textfont\bffam=\eighteenbf \scriptfont\bffam=\twelvebf
   \scriptscriptfont\bffam=\tenbf
\normalbaselineskip=20pt
  \setbox\strutbox=\hbox{\vrule height13.5pt depth6.5pt width0pt\relax}%
\everymath {\fam0 }
\everydisplay {\fam0 }
  \normalbaselines\bf}

\def\elevenbold{\def\rm{\fam0\elevenrm}%
  \textfont0=\elevenbf \scriptfont0=\eightbf \scriptscriptfont0=\sixbf%
  \textfont1=\elevenmib \scriptfont1=\eightmib \scriptscriptfont1=\sixmib%
  \textfont2=\elevensyb \scriptfont2=\eightsyb \scriptscriptfont2=\sixsyb%
  \textfont3=\elevenex \scriptfont3=\elevenex \scriptscriptfont3=\elevenex%
  \def\bf{\fam\bffam\elevenbf}%
  \textfont\bffam=\elevenbf \scriptfont\bffam=\eightbf%
   \scriptscriptfont\bffam=\sixbf%
\normalbaselineskip=14pt%
  \setbox\strutbox=\hbox{\vrule height10pt depth4pt width0pt\relax}%
\everymath {\fam0 }%
\everydisplay {\fam0 }%
  \normalbaselines\bf}


\hsize=125mm
\vsize=185mm

\parindent=8mm
\frenchspacing

\widowpenalty=10000
\clubpenalty=10000


\newif\ifappend

\def\appendix#1#2{\def\applett{#1}\def\two{#2}%
\global\appendtrue
\global\theoremcount=0
\global\eqcount=0
\vskip18pt plus 18pt
\vbox{\parindent=0pt
\everypar={\hskip\parfillskip}
\def\\ {\vskip1sp}\elevenpoint\bf Appendix%
\ifx\applett\empty\gdef\applett{A}\ifx\two\empty\else.\fi%
\else\ #1.\fi\hskip6pt#2\vskip12pt}%
\global\sectiontrue%
\everypar={\global\sectionfalse\everypar={}}\nobreak\ignorespaces}

\newif\ifRefsUsed
\long\def\references{\global\RefsUsedtrue\vskip21pt
\everypar={}
\vskip20pt\goodbreak\bgroup
\vbox{\centerline{\eightsc References}\vskip6pt}%
\ifdim\maxbibwidth>0pt\leftskip=\maxbibwidth\else\leftskip=18pt\fi%
\ifdim\maxbibwidth>0pt\parindent=-\maxbibwidth\else\parindent=-18pt\fi%
\everypar={\amref}%
\ninepoint
\frenchspacing
\nobreak\ignorespaces}

\def\endreferences{\vskip1sp\egroup\everypar={}%
\nobreak\vskip8pt\vbox{}}

\def\bibline{\hbox to30pt{\hrulefill}\/\/}

\def\name#1{{\eightsc#1}}

\newdimen\maxbibwidth
\def\AuthorRefNames [#1] {\def\amref{\spamref}
\setbox0=\hbox{[#1] }\global\maxbibwidth=\wd0\relax}

\def\spamref[#1] {\leavevmode\hbox to\maxbibwidth{\hss[#1]\hfill}}


\font\teneu=eufm10\font\seveneu=eufm7\font\fiveeu=eufm5
\font\tenlv=msbm11\font\sevenlv=msbm7\font\fivelv=msbm5

\newfam\eufam  \def\eu{\fam\eufam\teneu}
\textfont\eufam=\teneu \scriptfont\eufam=\seveneu 
                              \scriptscriptfont\eufam=\fiveeu

\newfam\lvfam  \def\lv{\fam\lvfam\tenlv}
\textfont\lvfam=\tenlv \scriptfont\lvfam=\sevenlv 
                              \scriptscriptfont\lvfam=\fivelv


\def\birmap{\hhb2\hbox{-\hhb{1}-\hhb{1}-}\hhb{-1}%
\lower-1pt\hbox{$\scriptscriptstyle>$}\hhb2}


\def\nix{{\phantom{|}}} 

\def\bsn{\bigskip\noindent}
 \def\msn{\medskip\noindent}
  \def\ssn{\smallskip\noindent}
   


  \def\hhb#1{\hbox to#1pt{}}
  \def\hor#1{\smash
          {\mathop{{\lgrghtar}}\limits^{\lower2pt\hbox{$\scriptstyle{#1}$}}}}
  \def\horr#1{\smash
         {\mathop{{\lglgrghtar}}\limits^{\lower2pt\hbox{$\scriptstyle{#1}$}}}}



\def\lgrghtar{{\hhb2{\relbar\joinrel\rightarrow}\hhb2}}
\def\lglgrghtar{{\hhb1{\relbar\joinrel\relbar\joinrel\rightarrow}\hhb1}}

\def\nmnm#1{{\elevensc #1}}


\def\plim#1{\hbox to14pt{lim\kern-13pt\lower4.5pt\hbox{$\scriptstyle
    \longleftarrow$}\kern-8pt\lower8.5pt\hbox{$\scriptstyle{#1}$}}}

\def\ilim#1{\hbox to14pt{lim\kern-13pt\lower4.5pt\hbox{$\scriptstyle
    \longrightarrow$}\kern-8pt\lower8.5pt\hbox{$\scriptstyle{#1}$}}}


\def\td{{\rm td}} 
\def\cd{{\rm cd}}  
\def\chr{{\rm char}}
\def\vid{{\lower-1pt\hbox{/}\kern-7pt{\hbox{O}}}}
\def\gata{} 

\def\ust#1{{}^*\!{#1}}
\def\gdw{if and only if\ }
\def\proof{{\it Proof}}

\def\grossl{\Lambda}
\def\kleinl{\lambda}
\def\eqv{ $\Leftrightarrow$ }


\def\chr{{\rm char}}
\def\cd{{\rm cd}} \def\vcd{{\rm vcd}}
\def\Th#1{{\eu Th}(#1)}
\def\abs#1{{#1}^{\rm abs}}
\def\wahr{{|\kern4pt\hbox{=}}}
\def\ee{elementarily equivalent}
\def\iff{if and only if\ }

\def\Hom{{\rm Hom}}
\def\rsdp{{\times\kern-3pt\lower-1pt\hbox{$\scriptscriptstyle|$}\,}}
\def\KM{{\rm K}^{\rm M}}
\def\pff#1{\;<\!<\!#1\!>\!>}

\def\pta#1{p^\nix_{a,#1}(T)}

\pageno=1
\vskip-20pt
\noindent
{\sixrm Variant of Nov 2000. Last reviews July/Sept 2001.} 
\bsn
\centerline{\eighteenbf Elementary equivalence}
\ssn
\centerline{\eighteenbf versus Isomorphism}

\bsn
\ssn
\centerline{ by Florian Pop at Bonn}


\bsn
\ssn
\centerline{\bf 1.\ Introduction}
\bsn
The aim of this note is to give some new evidence for a long standing 
open question concerning the relation between elementary equivalence 
and isomorphism in the class of finitely generated fields. More precisely,
we will give a {\it positive answer\/} to this question in the ``general''
case.
\ssn

We recall that two fields $K$ and $L$ are called \ee, if for every 
sentence $\psi$ in the language of fields on one has: $\psi$ is 
true in $K$ \iff $\psi$ is true in $L$. Just to start with, we recall
that the assertion ``$K$ admits an ordering'' can be expressed in an
elementary way (by a scheme of axioms) by saying that ``$-1$ is not 
a sum of squares in $K$'', wheres the fact that a function field 
$K|{\lv Q}$ has transcendence degree $d$ cannot be expressed in an 
obvious way (or not at all), as we would have to say that ``there 
exist $t_1,\dots,t_d\in K$ such that for all non-zero polynomials 
$p(X_1,\dots,X_d)$ over ${\lv Q}$ one has $p(t_1,\dots,t_d)\neq0$, and 
every other element of $K$ is algebraic over ${\lv Q}(t_1,\dots,t_d)$''. 
Clearly, this is not elementarily expressable --at least not-- in 
an obvious way. Concerning this we make the following remarks:

\ssn

{1)} For any field $K$ let $\Th K$ be the set of sentences in 
the language of fields, which are true in $K$. Then $K$ and $L$ are 
\ee\ if and only if $\Th K=\Th L$. Clearly, if $K\cong L$ as fields, 
then $\Th K=\Th L$, thus $K$ and $L$ are \ee. 

\ssn

{2)} The precise relation between elementary equivalence and 
isomorphism is as follows: $K$ and $L$ are \ee\ iff there exit 
ultra-powers $K^*=K^I/{\cal D}$ and $L^*=L^J/{\cal E}$ which are 
isomorphic as fields. In particular, if $K$ and $L$ are \ee, then 
they have the same prime field $k$, and if $\abs K$ and $\abs L$ 
denote their ``absolute subfields'' (i.e., the relative algebraic 
closures of $k$) in $K$ and in $L$ respectively, then $\abs K\cong 
\abs L$. 

\ssn 

{3)} From  a more geometric point of view, one can express the 
elementary equivalence as follows: $K$ and $L$ are \ee\ \iff they 
have the same prime field $k$, and for every constructible subset 
$S$ of a scheme of finite type over $k$ one has: $S(K)$ non-empty 
\iff $S(L)$ non-empty. 

\ssn

It appears to be an interesting and intriguing question,
whether {\it for finitely generated fields, the elementary
equivalence if the same as the isomorphism.\/} The first
non-obvious question in this direction is whether the 
absolute transcendence degree of a finitely generated 
field can be described by an assertion in the language of
fields. It seems that the following more precise question
was asked by \nmnm{Sabbagh} in the beginning if the 
Eighties: Let $K$ be a function field of one variable
over ${\lv Q}$, and $L$ the rational function field in two
variables over ${\lv Q}$. Is it then possible that $K$
and $L$ are \ee? Part of the problem here is clearly 
the question, whether one could find assertions in the 
language of fields which can distinguish between the 
transcendence degrees 1 and 2 of function fields over 
${\lv Q}$. Before giving the main results, let us recall 
that a function field $K|\kappa$ over some base field 
$\kappa$ is said to be of {\it general type,\/} if it is 
the function field of a projective smooth variety of general
type $X\to\kappa$. Correspondingly, we say that a finitely 
generated field $K$ is of general type, if it is of 
general type over its absolute subfield $\kappa=\abs K$. 
The main results of this note are the following:

\proclaim{Theorem A} ({\rm Arithmetic variant}). {\it Let $K$ 
and $L$ be finitely generated fields which are \ee. Let $\kappa$
and $\lambda$ be their absolute subfields. Then one has:
\ssn
\item{\rm(1)} $\kappa$ and $\lambda$ are isomorphic, and 
$\td(K|\kappa)$ equals $\td(L|\lambda)$.
\ssn
\item{\rm(2)} Moreover, there exits an embedding $\imath:K\to L$ 
such that $L$ is finite separable over $\imath(K)$. Furthermore, 
if $K$ is of general type, then $K\cong L$ as fields. 
\vskip0pt}

\proclaim{Theorem B} {(\rm Geometric variant}). {\it Let $K|\kappa$ 
and $L|\lambda$ be function fields over algebraically closed fields 
$\kappa$, respectively $\lambda$. Suppose that $K$ and $L$ are 
\ee\ as fields. The one has:
\ssn
\item{\rm(1)} $\kappa$ and $\lambda$ are \ee, and $\td(K|\kappa)$ 
equals $\td(L|\lambda)$. 
\ssn
\item{\rm(2)} Suppose $K|\kappa$ is of general type. Then there 
exist function subfields $K_0|\kappa_0\hookrightarrow K|\kappa$
and $L_0|\lambda_0\hookrightarrow L|\lambda$ such that $K=K_0\,\kappa$ 
and $L=L_0\,\lambda$, and $K_0|\kappa_0\cong L_0|\lambda_0$ as 
function fields.
\ssn
\item{} In particular, if $\kappa\cong\lambda$ are isomorphic, 
then $K|\kappa\cong L|\lambda$ are isomorphic as function fields. 
\ssn\/}

We would like to mention that there are already results 
in the literature concerning the geometric case by \nmnm{Duret} 
[D1], [D2], and \nmnm{Pierce} [Pi], which completely answer 
the case of function fields of curves, thus working under 
the hypothesis $\td(K|\kappa)=1=\td(L|\lambda)$. Among other 
things they show that the result is true also for the function 
fields of elliptic curves {\it without complex multiplication.\/} 
This goes beyond Theorem~B in the case under discussion, as 
the elliptic curves are not curves of general type.
 
\bigskip

It remains an {\bf open question} to give the precise 
{\it relation between elementary equivalence and isomorphism\/} 
in the ``non-general'' case.

\bigskip

We would also like to say that a preliminary form of the above
results, precisely: the characterisation of the transcendence
degree by Pfister forms was known to the author already in 1998;
and Theorem~A in the case $\td(K|\kappa)\leq1$, but the approach 
was different. The main tool was the Mordell Conjecture (as proved 
by \nmnm{Faltings}). In particular, we exploited the relation 
between {\it rational points\/} on general curves and the {\it 
elementary theory\/} of function fields of curves over number 
fields. It remains a serious {\bf open question} to precisely 
understand in general the relation between these two apparently 
different problems.

\bsn
\centerline{\bf Section 1. \ Detecting the transcendence degree}

\bsn
In this subsection we show how to detect the transcendence degree
-in the cases under discussion- by  ``nice'' sentences in the 
language of fields. In the positive characteristic case and in 
the geometric case, i.e., if the base field is algebraically
closed, the answer to this question is ``easy'', and might well
already be known. Therefore we will only indicate
how to get it. Nevertheless, in the arithmetic case, i.e., in 
the case of finitely generated fields over number fields, the 
only way we can do this is by using the the Milnor Conjecture, 
proved by \nmnm{Voevodsky, Rost, et al}, see \nmnm{Kahn}'s [Kh]
talk in S\'eminaire Bourbaki.

\smallskip

We consider the context: $K|\kappa$ is a function field
with $\td(K|\kappa)=d>0$. 
\msn
A) \ {\it $\kappa$ has positive characteristic\/}
\msn
Suppose $\kappa$ is a perfect field of characteristic $p>0$,
in particular $\kappa$ might be finite or algebraically closed.
Thus, if $K^p$ denotes the subfield of $K$ consisting of all the
$p$-powers in $K$, then $K|K^p$ is a (finite) purely inseparable
extension. Since $\kappa$ is a perfect field, it follows that
$[K:K^p]=p^d$, where $d$ is the transcendence degree of $K|\kappa$.
This fact can be easily expressed by a sentence in the language 
of fields. But a more precise formulation exists, by which we 
can say that a {\it given system\/} of elements of $K$ is a 
separable transcendence basis of $K|\kappa$ as follows. We 
recall the following two definitions: (1) First, a system 
of elements ${\cal B}_r=(t_1,\dots,t_r)$ of $K$ is said to be 
$p$-independent, if the system of all the monomials $\underline 
t^{\underline i}$ of the following
form $\underline t^{\underline i}=t_1^{i_1}\dots t_d^{i_r}$, 
$0\leq i_1,\dots,i_r<p$, is linearly independent over $K^p$.
(2) Second, a system of elements ${\cal B}_d=(t_1,\dots,t_d)$ 
of $K$ is said to be a $p$-basis, if the set of all the monomials 
as above is a vector space basis of $K$ over $K^p$. 
\par
It is clear that both assertion (1) and assertion (2) can be 
interpreted as formulas in the language of fields with no 
parameters excepting $t_1,\dots,t_r$. Namely for a given system 
${\cal B}_r=(t_1,\dots,t_r)$ of elements of $K$, let 
${\underline i}=(i_1,\dots,i_r)$ be all the multi-indexes with
$0\leq i_j<p$, and define the following form of degree $p$ in 
the $p^r$ variables $X_{\underline i}$ over $K$:
$$
q^{(p)}_{{\cal B}_r}(X_1,\dots,X_{p^r})=\sum_{\underline i}
       \underline t^{\underline i}X^p_{\underline i}
$$
Then the form $q^{(p)}_{{\cal B}_r}$ does not represent $0$ 
over $K$ if and only if ${\cal B}_r$ is $p$-independent; and 
the form $q^{(p)}_{{\cal B}_d}$ is universal over $K$, but 
does not represent $0$ over $K$ if and only if ${\cal B}_d$ is 
a $p$-basis of $K|\kappa$.
\ssn

The relation of this with the transcendence bases of $K|\kappa$
is given by the following well known fact, see e.g.\ [Ei], 
Appendix~1.
\msn
{\bf Fact 1.1.} \ {\it Let $K|\kappa$ be a function field with
$\kappa$ is a perfect field of characteristic $p>0$. For a system 
of elements ${\cal B}_r=(t_1,\dots,t_r)$ of $K$ the following 
assertions hold:
\ssn
\item{\rm(1)} ${\cal B}_r$ is $p$-independent \eqv ${\cal B}_r$
can be completed to a separable transcendence basis of $K|\kappa$
\eqv $q^{(p)}_{{\cal B}_r}$ does not represent $0$ over $K$.
\ssn
\item{\rm(2)} ${\cal B}_d$ is a $p$-basis \eqv ${\cal B}_d$ is
a separable transcendence basis of $K|\kappa$ \eqv the form
$q^{(p)}_{{\cal B}_d}$ is universal over $K$, but does not represent
$0$ over $K$.
\ssn 
\item{\rm(3)} In particular, $K|\kappa$ has transcendence degree $
d$ if and only if the following sentence is true in $K$: There exists 
${\cal B}_d=(t_1,\dots,t_d)$ such that $q^{(p)}_{{\cal B}_d}$ is 
universal over $K$, but does not represent $0$ over $K$. 
\ssn\/}

\vfill\eject
\noindent
B) \ {\it $\kappa$ is algebraically closed\/}
\msn
Let $K|\kappa$ be an a function field with $\kappa$ an 
algebraically closed field of characteristic $\neq2$. It
is interesting to remark that in the case $p=2$, the form
$q^\nix_{{\cal B}_r}:=q^{(2)}_{{\cal B}_r}$ defined above 
is a very special quadratic form, namely the $r$-fold 
Pfister form attached to the system of elements ${\cal B}_r$
of $K$. On the other hand, since $K|\kappa$ is a function 
field of transcendence degree $d$, and $\kappa$ is 
algebraically closed, $K$ is a $C_d$ field. For quadratic 
forms, this means that every such form in $n>2^d$ variables 
represents $0$ over $K$. In particular, if $\chr(\kappa)\neq2$,
every $d$-fold Pfister form $q=q^\nix_{\cal B}$ is universal. 
(In the case $\chr(\kappa)=2$ we know already what happens
from the previous subsection.) The point is that in this way 
we can characterize the transcendence degree of $K|\kappa$, 
and even more, determine whether a system of elements 
${\cal B}_d=(t_1,\dots,t_d)$ of $K$ is a transcendence basis
of $K|\kappa$.  

\msn
{\bf Fact 1.2.} \ {\it Let $K|\kappa$ be a function field
with $\kappa$ an algebraically closed field with $\chr(\kappa)
\neq2$. To systems ${\cal B}_r=(t_1,\dots,t_r)$ of elements $K$ 
we denote by $q^\nix_{{\cal B}_r}$ the corresponding Pfister 
form. Then the following holds:
\ssn
\item{\rm(1)} $K|\kappa$ has transcendence degree $d$ \eqv 
for all ${\cal B}_d$ the corresponding $q^\nix_{{\cal B}_d}$ 
is universal over $K$, and there exit ${\cal B}_d$ such that 
$q^\nix_{{\cal B}_d}$ does not represent $0$ over $K$.
\ssn
\item{\rm(2)} Suppose that $\td(K|\kappa)=d$. For a system of 
elements ${\cal B}_d$ of $K$ suppose that $q^\nix_{{\cal B}_d}$ 
is universal, but does not represent $0$ over $K$. Then 
${\cal B}_d$ is a transcendence basis of $K|\kappa$.
\ssn
\item{\rm(3)} In particular, let ${\cal S}\subset\kappa$ be an
infinite subset, and for a systems ${\cal B}_r$ as above, 
and $a=(a_1,\dots,a_r)$ set ${\cal B}_{r,a}=(t_1-a_1,
\dots,t_r-a_r)$. Then ${\cal B}_r$ is a transcendence base for 
$K|\kappa$ if and only if for ``almost all'' $a$ the corresponding 
form $q^\nix_{{\cal B}_{r,a}}$ is universal, but does not 
represent $0$ over $K$. 
\ssn\/}

{\it Proof.\/} 
First, since $\kappa$ is algebraically closed, $K$ is a $C_d$ 
field.  Thus, for every $a\in K^\times$, the quadratic form 
$q=aX^2-q^\nix_{{\cal B}_d}$ has $2^d+1$ variables, hence 
it represents $0$ over $K$. Since $\chr(\kappa)\neq2$, from 
this is follows that $q^\nix_{{\cal B}_d}$ represents $a$ over 
$K$. The remaining assertions can be derived as follows. Let
$X\to\kappa$ be some model of $K|\kappa$, and let $x\in X$ 
be a closed regular (thus smooth) point. Then for every regular 
system of parameters $(t_1,\dots,t_d)$ at $x$ one has the 
following

\ssn

{\it Claim:\/} \ The quadratic form $q^\nix_{(t_1,\dots,t_d)}$
does not represent $0$ over $K$.
\ssn

Indeed, it is a well known fact that, since $x$ is 
a regular closed point, $K$ has $\kappa$-embeddings in 
$\Lambda_d=\kappa((t_1))\dots((t_d))$. One shows now by induction 
on $r$, that the $r$-fold Pfister form $q^\nix_{(t_1,\dots,t_r)}
=\pff{t_1,\dots,t_r}$ is universal over $\Lambda_r$, but does 
not represent $0$. Finally, for the last assertion, remark 
that $r\leq{\rm td}(K|\kappa)$, as the form $q^\nix_{{\cal B}_{r,a}}$ 
does not represent $0$ over $K$. If now $r<{\rm td}(K|\kappa)$,
then for a proper choice of a closed point $x$ as above, after 
setting $a_i=t_i(x)$, it follows that $(t_1-a_1,\dots,t_r-a_r)$
is contained in a regular system of parameters at $x$. Thus if
$r<{\rm td}(K|\kappa)$, and $(t_{r+1},\dots,t_d)$ are the
extra-parameters, then none of these is represented by 
$q^\nix_{{\cal B}_{r,a}}$ over $\Lambda_{d,a}=\kappa((t_1-a_1,))
\dots((t_d))$, thus nor over $K$. 
\gata 
\msn
C) \ {\it $\kappa$ is a number field\/}

\ssn
We now come to the discussion of the arithmetic case, which
is the most interesting (but the more difficult) one. Thus here
$K|\kappa$ is a function field over a number field $\kappa$.
As above, to systems ${\cal B}_r=(t_1,\dots,t_r)$ of elements 
of $K$, we denote by $q^\nix_{{\cal B}_r}$ the corresponding 
Pfister form. 
\msn
{\bf Fact 1.3.} \ {\it Let $K|\kappa$ be a function field
with $\kappa$ a number field. Then the following holds:
\ssn
\item{\rm(1)} $K|\kappa$ has transcendence degree $d$ \eqv \
for all ${\cal B}_{d+2}$ the corresponding quadratic form 
$q^\nix_{{\cal B}_{d+2}}$ is universal over $K[\sqrt{-1}]$,
and there exits ${\cal B}_{d+2}$ such that $q^\nix_{{\cal B}
_{d+1}}$ does not represent $0$ over $K[\sqrt{-1}]$.
\ssn
\item{\rm(2)} Suppose that $\td(K|\kappa)=d$. For a system 
of elements ${\cal B}_d=(t_1,\dots,t_d)$ of $K$ suppose that 
there exit $t_{d+1},t_{d+2}\in\kappa$ such that $q^\nix_{(t_1,
\dots,t_{d+2})}$ does not represent $0$ over $K[\sqrt{-1}]$. 
Then ${\cal B}_d$ is a transcendence basis of $K|\kappa$.
\ssn\/}

\proof. \ We proof goes along the following steps.
\ssn
a) \ {\it Some cohomological computations\/}
\msn 
For a field extension $K|\kappa$ as above, and a prime number 
$\ell$, we denote by $\vcd_\ell(K)$ the virtual $\ell$-cohomological 
dimension of $K$. It actually equals $\cd_\ell(E)$ for 
every non-real finite extension $E|K$. In particular, if 
$\sqrt{-1}$ lies in the field in discussion, then $\vcd_\ell=
\cd_\ell$ for every $\ell$. 
\ssn

By classical results in Galois cohomology of fields, see
e.g.\ \nmnm{Serre} [S], Ch.II, one has: First, if $\ell=\chr(K)$, 
then $\cd_\ell(K)=1$ (this is a theorem of Shafarevich). Second, 
if $K$ is finite, then $\cd_\ell(K)=1$ for all $\ell$, and if $K$ 
is a number field, then $\vcd_\ell(K)=2$ for all $\ell$ (this is
a theorem of Tate). Combined with the fact that $\vcd_\ell
\big(\kappa(t)\big)=\vcd_\ell(\kappa)+1$ for all fields $\kappa$ 
and $\ell\neq\chr(\kappa)$, we get:

\msn
{\bf Fact.\/} \ {\it In the context of Fact~1.3 one has: 
$\td(K|\kappa)=d=\vcd_\ell(K)-2$ for all rational prime 
numbers $\ell$.\/}
\msn

Let $\mu_\ell$ be the group of $\ell^{\rm th}$ roots of unity.
Then by Kummer theory and the tame symbol from the Milnor 
$K$-theory, we have a canonical isomorphism 
$
\delta:K^\times/\ell\cong{\rm H}^1(G_K,\mu_\ell)
$, 
which induces canonical homomorphisms
$$
h:K^\times\otimes\dots\otimes K^\times\to
        \KM_n(K)/\ell\hor{d}{\rm H}^n(G_K,\mu_\ell^{\otimes n}),
$$ 
denoted 
$
(a_1,\dots,a_n)\mapsto\{a_1,\dots,a_n\}\mapsto a_1\cup\dots\cup a_n,
$
see e.g.\ \nmnm{Milnor} [M1]. We further remark, that if 
$\mu_\ell\subset K$, then 
${\rm H}^1(G_K,\mu_\ell)\cong\Hom(G_K,{\lv Z}/\ell)$ is the 
character group of $G_K$ with values in ${\lv Z}/\ell$, and 
second, ${\rm H}^d(G_K,\mu_\ell^{\otimes n})$ can be identified with 
${\rm H}^n(G_K,{\lv Z}/\ell)$. The main observation is now 
the following 

\msn
{\bf Lemma}. {\it In the context of Fact~1.3, let $\ell$ be 
a prime number. Let $E$ be a finite extension of $K$ containing 
$\mu_{2\ell}$. Then one can detect $m=\vcd_\ell(K)$ as the 
unique natural number $m$ with the following properties:
\ssn
\item{\rm(i)} There exist $a_1,\dots,a_m$ in $K^\times$ 
such that $a_1\cup\dots\cup a_m\neq0$ as element of 
${\rm H}^m(G_E,{\lv Z}/\ell)$.  
\ssn
\item{\rm(ii)} If $a_1\dots,a_{m+1}$ are arbitrary elements of 
$K^\times$, then $a_1\cup\dots\cup a_{m+1}=0$ in  
${\rm H}^{m+1}(G_E,{\lv Z}/\ell)$.
\ssn\/}

{\it Proof.\/} Under the hypothesis on the Lemma on has 
$\vcd_\ell(E)=\cd_\ell(E)$. By the discussion preceeding the Lemma, 
the only less obvious fact is to show that for $m=d+2$ the 
assertion (i) is satisfied. We show that this is true as follows: 
Let $Y$ be a model of $K$ over ${\cal O}_\kappa$. Let $x\in Y$ 
be a smooth point of $Y\to{\cal O}_\kappa$ with the properties: 

\item{(i)} Its residue field $\kappa^\nix_x$ satisfies 
$\chr(\kappa^\nix_x)\neq\ell$ and contains $\mu_{2\ell}$. 

\item{(ii)} $x$ is totally split in the field extension $E|K$.

\noindent
(Such points $x$ do exist, why?) Let $(t_{d+1},t_d,\dots,t_1)$ 
be a regular sequence at $x$ such that $t_{d+1}\in{\cal O}_\kappa$. 
(This is possible, as $Y\to{\cal O}_\kappa$ is smooth at $x$.) 
The residue field $\kappa^\nix_x$ is a finite field. Choose 
$\overline t_{d+2}\in\kappa^\nix_x$ a non-$\ell^{\rm th}$ power, 
and $t_{d+2}\in K$ a representative for it. 
\msn

{\it Claim:\/} \ $t_1\cup\dots\cup t_{d+2}\neq0$ \ in \ 
   ${\rm H}^{d+2}(G_E,{\lv Z}/\ell)$.
\msn
The idea of the proof is to construct an algebraic extension
$\Lambda$ of $K$ containing $E$, and show that the restriction 
of $t_1\cup\dots\cup t_{d+2}$ to $\Lambda$ is $\neq0$, thus the same 
must hold over $E$. We take for $\Lambda$ a ``multi-local
field'' at $x$. If namely ${\eu p}$ is the place of $\kappa$ 
determined by the image of $x$ in ${\rm Spec}\,{\cal O}_\kappa$,
and $\kappa^\nix_{\eu p}$ is the completion of $\kappa$ at ${\eu p}$,
then $\Lambda\cong_\kappa \kappa^\nix_{\eu p}((t_d))\dots((t_1))$.
Thus, if $v$ is the valuation of $\Lambda$, then 
$$
v(K)=v(\Lambda)=:\Gamma\cong{\lv Z}^{d+1}
$$ 
with generators $v(t_{d+1}), v(t_d),\dots,v(t_1)$. We also 
remark that the relative algebraic closure $K^{\rm h}$ of $K$ 
in $\Lambda$ is a Henselisation of $K$ at $v$. The absolute 
Galois group of $K^{\rm h}$ is canonically isomorphic to the 
absolute Galois group $Z=G_\Lambda$ of $\Lambda$. By the 
hypothesis on $x$, it follows that $E$ is $K$-embedable in 
$K^{\rm h}\subset\Lambda$. Further, by general valuation 
theory, see e.g.\ the explanations in [P], Local Theory, 
since $\mu_{2\ell}\subset\Lambda$, we have:
\ssn

{1)} \ Let $Z(\ell)$ and $G(\ell)$ be the maximal pro-$\ell$ quotients 
of $Z$, respectively of $G_{\kappa_x}\cong\widehat{\lv Z}$. 
Then $Z(\ell)$ has a well known structure given as follows: 
If $T(\ell)$ denotes the inertia subgroup of $Z(\ell)$, then 
one has a split group extension
$$
1\to T(\ell)\to Z(\ell)\to G(\ell)\to1.
$$
Further, $T(\ell)$ is canonically ismorphic as a $G(\ell$)-module 
to ${\rm Hom}(\widetilde\Gamma,\mu_\infty)$, where 
$\widetilde\Gamma=(\Gamma\!\otimes\!{\lv Q})\,/\,\Gamma$, and 
$\mu_\infty$ is the group of roots of unity whose orders are 
$\ell$-powers. Thus $Z(\ell)$ is isomorphic as an abstract 
profinite group to ${\lv Z}_\ell^{d+1}\rsdp G(\ell)$, where 
$G(\ell)\cong{\lv Z}_\ell$ acts on ${\lv Z}_\ell^{d+1}$ 
componentwise via the cyclotomic character of $G(\ell)$.
\ssn

{2)} \ Let $Z_\ell$ denote a Sylow $\ell$-group of $Z$. Then the 
canonical projection $Z_\ell\to Z(\ell)$ is an isomorphism. See
loc.cit..
\ssn

In particular, the isomorphy type of $Z(\ell)$ depends only 
on $\vcd_\ell(K)$, and on the (cyclotomic character of the) 
residue field $\kappa^\nix_x$ of $x$. For the given residue field 
$\kappa^\nix_x$, thus the cyclotomic character of $G(\ell)$, 
and a given finite free ${\lv Z}_\ell$-module $\Delta_n$ of 
rang $(n-1)$ on which $G(\ell)$ acts via the cyclotomic character, 
we set $Z_n=\Delta_n\rsdp G(\ell)$. Thus the above $Z(\ell)$ is
isomorphic as a profinite group to $Z_{d+2}$. The proof of the
Claim will then follow from the following group theoretical 
fact:

\msn
{\bf Fact.} \ {\it One has $\cd_\ell(Z_n)=n$. More precisely,
if $0\neq\chi_n\in{\rm H}^1(G(\ell),{\lv Z}/\ell)$, and further
$\chi_1,\dots,\chi_{n-1}\in{\rm H}^1(\Delta_{n},{\lv Z}/\ell)$ 
are independent characters, then $\chi_1\cup\dots\cup\chi_n\neq0$ 
and generates ${\rm H}^n(Z_n,{\lv Z}/\ell)$.\/}
\msn

This is a well known fact, thus only sketch here a proof: 
Consider an exact sequence of the form 
$
0\to\Delta\to\Delta_n\to\Delta_{n-1}\to0$ such that 
$\Delta\subset\ker(\chi_i)$ for $i<n-1$, and in 
particular, $\chi_{n-1}$ non-zero on $\Delta$. The 
above exact sequence gives rise to an exact sequence 
$0\to\Delta\to Z_n\to Z_{n-1}\to 1$. Thus using for 
example a Hochschild--Serre spectral sequences of the form
${\rm H}^p\big(Z_{n-1},{\rm H}^q\big(\Delta,{\lv Z}/\ell)\big)
\Rightarrow H^{p+q}\big(Z_n,{\lv Z}/\ell\big)$, with a trivial
action of $\Delta$ on ${\lv Z}/\ell$, one immediately obtains 
by induction the following: ${\rm H}^n\big(Z_n,{\lv Z}/\ell)$
is non-trivial, and it is generated by the cup-product 
$\chi_1\cup\chi_2\cup\dots\cup\chi_n$.
\ssn

To finish the proof of the above claim, we remark that in our 
context we can coose $\chi_i=\delta(t_i)$, $i=1,\dots,d+2$, 
where $\delta$ is the Kummer isomorphism. The hypotheses on 
the $\chi_i$ are satisfied, as by thier choice the elements 
$t_1,\dots,t_{d+2}$ are independent in $\Lambda^\times/\ell$, and
moreover, build an ${\lv F}_\ell$-basis of $\Lambda^\times/\ell$. 
Finally, point 2) above implies that ${\rm H}^i(Z,{\lv Z}/\ell)$
and ${\rm H}^i\big(Z(\ell),{\lv Z}/\ell\big)$ are isomorphic.

\msn
b) {\it Transcendence degree and arithmetic\/}
\msn
As mentioned above, for every field $E$ there exists a canonical 
homomorphism called the tame residue symbol: 
$$
h_{\ell,n}:\KM_n(E)\,/\,\ell\to {\rm H}^n(G_K,\mu_\ell^{\otimes n})
$$
It is conjectured that $h_{\ell,n}$ is an isomorphism for all $\ell$ 
prime to $\chr(E)$. This is a generalisation of the so called 
{\it Milnor Conjecture,\/} which is the above assertion for $\ell=2$. 
The point is that the Milnor Conjecture has
a deep arithmetic significance relatetd to the arithmetic of the
quadratic forms (for general $\ell$ we do not have yet an
interpretation of the ``generalised Minor Conjecture''). The 
Milnor conjecture is now proved by contributions of several people,
with last major steps being done by \nmnm{Voevodsky, Rost, et
al}, see \nmnm{Kahn} [Kh]. We describe below the facts which 
are significant for us. 
\ssn

With $E$ as above, let $W(E)$ be the Witt ring of $E$, i.e., the
set of the isomorphy classes of anisotropic quadratic froms over
$E$ with the usual addition and multiplication. Let $I(E)$ be 
the ideal of even dimensional quadradic forms, and $I^n(E)$ its
powers. For $a_1,\dots,a_n$ in $E^\times$ let $\pff{a_1,\dots,a_n}$ 
be the corresponding ($n$-fold) Pfister form. The set of all
$n$-fold Pfister forms generates $I^n(E)/I^{n+1}(E)$. Milnor defined 
for every $n$ a homomorphism $d_n:\KM_n(E)/2\to I^n(E)/I^{n+1}(E)$ 
and conjectured that both $d_n$ and $h_n:=h_{2,n}$ are isomorphisms. 
In particular, this would give rise to {\it higher cohomological 
invariants\/} for quadratic forms, as we would then have then 
isomorphism $e_n:I^n(E)/I^{n+1}(E)\to{\rm H}^n(G_K,{\lv Z}/2)$ 
for every $n$, thus generalising the ${\rm dim}({\rm mod}\,2)$, the 
discriminant, and the Clifford (thus Hasse-Witt) invariant. 
The above isomorphism would work at the level of the Pfister 
forms as follows: 
$$
e_n:\;\;\pff{a_1,\dots,a_n}\;\>\mapsto (-a_1)\cup\dots\cup(-a_n).
$$   
(Note that the minus-sign comes from a convention which is not
necessarily the same in all sources. It depends on the definition
of $\pff{a}$ which is $\pff{a}\;=X_0^2\pm aX_1^2$. We work with the
``+'' convention.) See \nmnm{Elman--Lam} [E--L], \nmnm{Jacobs--Rost}, 
[J--R], etc.\ for more details and literature.
\ssn

We recall the following fact: Let $q=\pff{a_1,\dots,a_n}$ be a 
given Pfister form, and $a\in E^\times$. Then $q$ represents $-a$ 
\iff $q\,\otimes\hbox{$\pff{a}$}$ is hyperbolic. Thus passing to 
Galois cohomology and using the Milnor Conjecture we get the following:
\msn
{\bf Fact.} \ {\it For a Pfister form $q=\pff{a_1,\dots,a_n}$ and 
$a\in E^\times$ the following are equivalent:
\ssn
\item{\rm(1)} $q$ represents $-a$.
\ssn
\item{\rm(2)} $q\;\otimes\!\pff{a}\;=\pff{a_1,\dots,a_n,a}$ 
is hyperbolic.
\ssn
\item{\rm(3)} $(-a_1)\cup\dots\cup(-a_n)\cup(-a)=0$ 
in ${\rm H}^{n+1}(G_K,{\lv Z}/2)$.
\ssn
\/}

\par

Coming back to the proof of Fact~1.3, we remark that assertion~(1)
immediately follows from the discussion above with $E=K[\sqrt{-1}]$.
To prove assertion~(2), one proceeds by contradiction: Suppose 
that ${\cal B}_d=(t_1,\dots,t_d)$ is not a transcendence basis 
of $K|\kappa$, and let $K_0=\kappa({\cal B})$ be the subfield of 
$K$ generated by ${\cal B}$ over $\kappa$. Then $t(K_0|\kappa)<d$.
Thus for all $t_{d+1},t_{d+2}\in\kappa$ on has: The Pfister form
$q^\nix_{(t_1,\dots,t_{d+1})}$ is universal over $K_0[\sqrt{-1}]$,
in particular, it represents $t_{d+2}$. Therefore, $q^\nix_
{(t_1,\dots,t_{d+2})}$ is hyperbolic over $K_0[\sqrt{-1}]$,
thus over $K[\sqrt{-1}]$, contradiction~!
\gata

\msn

Summarizing, we have the following way to describe the transcendence 
degree and even more, transcendence bases of a function field.

\ssn
\proclaim{Theorem 1.4}. \ {\it Let $K|\kappa$ be a regular 
function field over some base field $\kappa$, where $\kappa$ 
is either algebraically closed, or finite, or a number field. 
To a system ${\cal B}=(t_1,\dots,t_r)$ of elements of $K$ and 
a rational prime number $p$ let 
$$
q^{(p)}_{{\cal B}_r}={\textstyle\sum_{\underline i}}\,
    {\underline t}^{\underline i}\,X_{\underline i}^p
$$ 
be the homogenous form over $K$ whose coefficients are all
the monomials of the form ${\underline t}^{\underline i}=
t_1^{i_1}\dots t_r^{i_r}$ with $0\leq i_1,\dots,i_r<p$.
Thus $q^\nix_{{\cal B}_r}:=q^{(2)}_{{\cal B}_r}$ is the 
$r$-fold Pfister form defined by ${\cal B}$. Then one has:
\ssn
\item{\rm(1)} Suppose $\kappa$ has characteristic $\chr(\kappa)=p$. 
Then ${\cal B}_d$ is a separable transcendence basis of $K|\kappa$
\gdw $q^{(p)}_{{\cal B}_d}$ is universal, but does not represent 
$0$ over $K$.
\ssn
\item{\rm(2)} Suppose $\kappa$ is algebraically closed with 
$\chr(\kappa)\neq2$. Then $\td(K|\kappa)=d$ \gdw $q^\nix_{{\cal B}_d}$ 
is universal over $K$ for every ${\cal B}_d$, and there does exist
${\cal B}_d$ such that $q^\nix_{{\cal B}_d}$ does not represent 
$0$ over $K$.
\ssn
\item{\rm(3)} Suppose $\kappa$ is number field. Then 
$\td(K|\kappa)=d$ \gdw $q^\nix_{{\cal B}_{d+2}}$ is universal 
over $K[\sqrt{-1}]$ for every ${\cal B}_{d+2}$, and there does
exist ${\cal B}_{d+2}$ such that $q^\nix_{{\cal B}_{d+2}}$ 
does not represent $0$ over $K[\sqrt{-1}]$.
\ssn
\item{} Moreover, $(t_1,\dots,t_d)$ is a transcendence basis
of $K|\kappa$ if there exit $t_{d+1},t_{d+2}\in\kappa$ such
that $q^\nix_{(t_1,\dots,t_{d+2})}$ does not represent $0$
over $K[\sqrt{-1}]$.
\ssn
In particular, if $L|\lambda$ is another function field with
$\lambda$ as above, such that $K$ and $L$ are \ee, then one 
has: $\kappa$ and $\lambda$ are \ee, thus they have the same 
characteristic; and they are isomorphic if $\kappa$ is finite 
or a number field; further $\td(K|\kappa)=\td(L|\lambda)$. \/}
\par

\proof. \ The assertions (1), (2), (3), are nothing but 
Facts 1.1, 1.2, 1.3 above. Concerning the last assertions, 
if $K$ and $L$ are \ee, then $\abs\kappa=\abs K\cong
\abs L=\abs\lambda$.

\bsn
\centerline{\bf 2.\ Generalities about function fields}

\bsn
In this section we put together well known facts about function
fields and their models, facts which will be used later.
\msn
A) \ {\it Rings of definition\/}
\msn
Let $K|\kappa$ be an arbitrary, separably generated function 
field with $\kappa$ the constant field of $K$. Further let
$X={\rm Proj}\Big(\kappa[\underline X]/(\underline f)\Big)$ 
be some projective model of $K|\kappa$. Here $\underline 
X=(X_0,\dots,X_n)$ is a system of $(n+1)$ indeterminates, and 
$\underline f=(f_1,\dots,f_m)$ is a system of homogenous 
polynomials in $\underline X$. We denote $\underline x=(x_0,
\dots,x_n)$ the resulting homogeneous coordinates on $X$.

\ssn

In the above context, let $k$ be the prime field of $\kappa$.
We define a {\it ring of definition\/} of $X\to\kappa$, and 
thus of $K|\kappa$, to be a $k$-subalgebra of finite type 
$R\subset\kappa$ over which the variety $X\to\kappa$ ``is 
defined''. The first approximation in doing this, is to take 
a $k$-subalgebra of finite type $R$ of $\kappa$ over which 
the ideal $(\,\underline f\,)$ is defined, and set ${\cal X}=
{\rm Proj}\Big(R[\underline X]/(\,\underline f\,)\Big)$. 
By replacing $R$ by a bigger $k$-subalgebra of finite type of 
$\kappa$, we can suppose that ${\cal X}\to R$ and its fibers
$X_s={\cal X}\times_R\kappa^\nix_s\to\kappa^\nix_s$ at points 
$s\in{\rm Spec}(R)$ have as nice properties as the generic fiber 
$X\to\kappa$ has. Thus we can/will suppose that:
\ssn
\item{0)} ${\cal X}\to R$ is flat. 
\ssn
\item{1)} Further, if $\kappa^\nix_s$ is the residue field at 
a point $s$ in ${\rm Spec}(R)$, then the fiber $X_s={\rm Proj}
\Big(\kappa^\nix_s[\underline X]/(\underline f_s)\Big)\to
\kappa^\nix_s$ is geometrically irreducible. 
\ssn
Let ${\cal B}=\underline t=(t_1,\dots,t_d)$ be a fixed 
(separable) transcendence basis of $K|\kappa$, and set 
$t_\nu=g_\nu/h_\nu$ with $g_\nu$ and $h_\nu$ homogeneous 
elements of $\kappa[\underline X]/(\underline f)$ of 
the same degree. Then, after enlarging $R$, we can suppose 
that all the $g_\nu$ and $h_\nu$ are defined over $R$. We 
will say in this case that ${\cal B}$ is defined over $R$. 
Further, if $g_{\nu,s}$ and $h_{\nu,s}$ are the specializations 
of $g_\nu$ respectively $h_\nu$ at some $s\in{\rm Spec}(R)$,
and the resulting rational functions $t_{\nu,s}=g_{\nu,s}/h_{\nu,s}$ 
on $X_s$ are defined, then we can consider the system 
${\cal B}_s=\underline t_s=(t_{1,s},\dots,t_{d,s})$ of 
rational functions in in the function field $K_s=\kappa(X_s)$ 
over $\kappa^\nix_s$. We will say that in this case that 
${\cal B}_s$ is defined at $s$. In this context, after 
shrinking ${\rm Spec}\,R$, we can suppose that for all 
$s\in{\rm Spec}(R)$ it holds: 
\ssn
\item{2)}  If $t_{\nu,s}=g_{\nu,s}/h_{\nu,s}$ are defined, and 
the system ${\cal B}_s=\underline t_s=(t_{1,s},\dots,t_{d,s})$ 
is a (separable) transcendence basis of the function field 
$K_s=\kappa^\nix_s(X_s)$ over $\kappa^\nix_s$.
\ssn
\item{} In particular, we can/will suppose that the above 
condition is satisfied at every $s\in{\rm Spec}(R)$, and that 
we have $[K:\kappa({\cal B})]=[K_s:\kappa^\nix_s({\cal B}_s)]$.    
\ssn
\item{3)} If $X$ is (projectively) normal, then the same is
true for ${\cal X}$, and for its fibers $X_s$.
\ssn
\item{4)} If $X\to\kappa$ is smooth, the same is true for 
${\cal X}\to R$. In particular, the fibers $X_s\to\kappa^\nix_s$ 
are then projective and smooth.  
\msn
Now suppose that $X\to\kappa$ is smooth, and by
the remarks above, the same is true for ${\cal X}\to R$.
This means that the sheaf of the relative differentials
$\Omega_{{\cal X}|R}$ is locally free of rank $d=\td(K|\kappa)$,
thus its $d^{\rm th}$ exterior power $\omega_{{\cal X}|R}$
is an invertible sheaf on ${\cal X}$ (the relative canonical
class). Taking this into account and using the Semi-continuity 
Theorem, after again shrinking ${\rm Spec}\,R$, we get: 
\ssn
\item{5)} Given a finite set ${\eu I}$ of birational invariants, 
like for instance plurigenera and/or Hodge numbers of $X$, it follows 
that the special fibers $X_s$ have the same birational invariants.  
\ssn
\item{6)} In particular, if ${\cal D}$ is a (relative) effective
divisor in the divisor class of $\omega_{{\cal X}|R}$, and $D$
and $D_s$ are its generic fiber respectively its fiber at $s$,
then the linear spaces $|m{\cal D}|$, $|mD|$, and $|mD_s|$
have the same dimension (for $m$ in some range $0\leq m\leq n$),
and the canonical $R$-morphism ${\cal X}\to{\lv P}^N_R$
with $N$ the dimension of $|m{\cal D}|$, give rise by base change
to the corresponding canonical projective $\kappa$-morphism
$X\to{\lv P}^N_\kappa$ respectively $\kappa^\nix_s$-morphism
$X_s\to{\lv P}^N_{\kappa_s}$.

\msn
B) \ {\it Approximations of function fields\/} 
\msn
Let $K|\kappa$ and $\grossl|\kleinl$ be function fields. 

\ssn

Let ${\cal X}\to R$ be a model of $K|\kappa$ as above.
We define an {\it approximation of $K|\kappa$ with 
values in $\grossl|\kleinl$ via ${\cal X}\to R$\/} to 
be a $\grossl$-rational point  
$
\varphi:{\rm Spec}\,\grossl\to{\cal X}
$ 
of ${\cal X}$ with the property: There exists a morphism 
of rings $p:R\to\kleinl$ such that denoting 
$X_\varphi={\cal X}\times_R\kleinl$ the base change of 
${\cal X}$ to $\kleinl$ via $p$, the point $\varphi$ is 
defined by a $\kleinl$-embedding of function fields 
$\imath^\nix_\varphi:\kleinl(X_{\kleinl})\to\grossl$.
This implies in particular that 
$
\td(K|\kappa)=\td(\kleinl(X_\varphi)|\kleinl)
                 \leq\td(\grossl|\kleinl)
$. 
For an approximation $\varphi$ of $K|\kappa$ with values 
in $\grossl|\kleinl$ via some given model ${\cal X}\to R$ 
as above we now make the following definitions:
 
\ssn
\item{-} Let ${\cal B}_d=(t_1,\dots,t_d)$ be a separable 
transcendence base of $K|\kappa$ as in subsection A), 2) 
above. We say that $\varphi$ is ${\cal B}$-separable, if 
$\grossl|\lambda$ is separably generated, and
$\imath^\nix_\varphi({\cal B})$ is contained in a separable
transcendence base of $\grossl|\kleinl$.

\ssn
\item{-} We will say that $\varphi$ is smooth, if ${\cal X}\to R$ 
is a smooth morphism. In particular, $X_\varphi\to\kleinl$ is
a smooth morphism.

\ssn
\item{-} We will say that $\varphi$ is minimal, if $\dim(R)$ 
is minimal among the dimensions of all the possible rings
of definition of $K|\kappa$. 

\ssn 
\item{-} We will say that $K|\kappa$ has {\it enough
(separable) approximations\/} with values in $\grossl|\kleinl$, 
if for every model ${\cal X}\to R$ of $K|\kappa$ (endowed 
with a separable transcendence basis ${\cal B}$) as above, 
$K|\kappa$ has (${\cal B}$-separable) approximations with 
values in $\grossl|\kleinl$ via ${\cal X}\to R$.

\msn

\proclaim{Theorem 2.1}. \ {\it In the above notations
the following holds:
\ssn
\item{\rm(1)} Let $K|\kappa$, $\grossl|\kleinl$, 
$\grossl'|\kleinl'$ be function fields. Suppose that 
$K|\kappa$ has enough (separable) approximations with 
values in $\grossl|\kleinl$, and that $\grossl|\kleinl$  
has enough (separable) approximations with values in 
$\grossl'|\kleinl'$. Then $K|\kappa$ has enough (separable)
approximations with values in $\grossl'|\kleinl'$.
\ssn
\item{\rm(2)} Let $K|\kappa$ and $L|\lambda$ be function 
fields such that both: $K|\kappa$ has enough separable 
smooth approximations with values in $L|\lambda$, and 
$L|\lambda$ has enough separable approximations with 
values in $K|\kappa$. Further suppose that $K|\kappa$ 
is a function field of general type. Then $K|\kappa$ and 
$L|\lambda$ have isomorphic function subfields $\imath_0:
K_0|\kappa_0\to L_0|\lambda_0$ such that $K=K_0\,\kappa$ and 
$L=L_0\,\lambda$.
\ssn
\item{} In particular, if the isomorphism $\imath_0:\kappa_0
\to\lambda_0$ can be prolonged to an isomorphism $\kappa\cong
\lambda$, then $\imath_0$ has a prolongation to an isomorphism
of function fields $\imath:K|\kappa\to L|\lambda$.
\ssn
\/}

\proof. \ To (1): Let $R\subset\kappa$, ${\cal X}\to R$ 
(and a transcendence base ${\cal B}$ of $K|\kappa$ as in
subsection A), 2) above) be given. We show that $K|\kappa$ 
has (${\cal B}$-separable) approximations with values 
in $\grossl'|\kleinl'$ via ${\cal X}\to R$. In the above
notations, let $\varphi:\grossl\to{\cal X}$ be a 
(${\cal B}$-separable) approximation of $K|\kappa$ with 
values in $\grossl|\kleinl$. Let $p:R\to\kleinl$ be a 
ring homomorphism defined by $\varphi$, and 
$\imath^\nix_\varphi:\kleinl(X_\varphi)\to\grossl|\kleinl$
the $\kleinl$-embedding of function fields defined by 
$\varphi$. By general facts concerning rational dominant 
maps we get: Given any model $Y\to\kleinl$ of $\grossl|\kleinl$, 
the morphism of function fields $\imath^\nix_\varphi$ is 
defined by a rational map $\phi:Y\birmap X_\varphi$. 
Moreover, there exist projective models $Y\to\kleinl$ 
of $\grossl|\kleinl$ such that $\imath^\nix_\varphi$ 
is induced by a dominant $\kleinl$-morphism
$\phi^\nix_\varphi:Y\to X_\varphi$. Let now $S\subset\kleinl$ 
be a finitely generated $k$-subalgebra containing $p(R)$, 
over which $Y$ and the morphism $\phi^\nix_\varphi$ are 
defined. This means, there exits a model ${\cal Y}\to S$ 
of $\grossl|\kleinl$ whose base change to $\kleinl$ is 
$Y\to\kleinl$, togehter with a morphism $\phi:{\cal Y}
\to{\cal X}$ which defines $\varphi$ as the composition
$$
\varphi:{\rm Spec}\,\grossl\hhb3\horr{{\rm can}}\hhb3 
       {\cal Y}\hhb3\horr{\phi}\hhb3{\cal X}
$$ 
We can/will suppose that ${\cal Y}\to S$ satisfies the 
conditions 0) and 1) from subsection A) above. Further, 
if the approximation $\varphi$ was ${\cal B}$-separable, 
let ${\cal C}$ be a separable transcendence basis of 
$\grossl|\kleinl$ containing $\imath^\nix_\varphi({\cal B})$. 
Then we can/will suppose that condition 2) of subsection A) 
is satisfied for ${\cal Y}\to S$ endowed with ${\cal C}$. 

\smallskip

Now let ${\cal X}_{S}={\cal X}\times_RS$ be the base 
change of ${\cal X}\to R$ to $S$ via $p:R\to S$. Then 
there exists a dominant $S$-morphism $\phi^\nix_S:
{\cal Y}\to{\cal X}_S$, whose base change to $\kleinl$ is 
the given $\kleinl$-morphism $\phi^\nix_\varphi:Y\to X_\varphi$. 
Thus in particular, $\phi^\nix_S$ defines the function field 
homomorphism $\imath^\nix_\varphi$.

\par

In order to finish the proof, let $\varpi:{\rm Spec}\,
\grossl'\to{\cal Y}$ be a (${\cal C}$-separable) 
approximation of $\grossl|\kleinl$ with values in 
$\grossl'|\kleinl'$ via ${\cal Y}\to S$. On then 
checks without difficulty that the composition
$$
\varphi':{\rm Spec}\,\grossl'\hhb5\horr{\varpi}
       \hhb5{\cal Y}\hhb5\horr{\phi_S}\hhb5{\cal X}_S
                 \hhb5\horr{{\rm can}}\hhb5{\cal X}
$$
is a $({\cal B}$-separable) approximation of $K|\kappa$ 
with values in $\grossl'|\kleinl'$ via ${\cal X}\to R$. 

\ssn
To (2). We first remark that the existence of
approximations of both: $K|\kappa$ with values in $L|\lambda$,
and vice-versa, implies $\td(K|\kappa)=\td(L|\lambda)$.

\par

In the notations from the proof of (1) above, we replace 
$\grossl|\kleinl$ by $L|\lambda$, and $\grossl'|\kleinl'$ 
by $K|\kappa$. For all possible separable transcendence bases 
${\cal B}$ of $K|\kappa$, we consider all ${\cal B}$-separable 
approximations of $K|\kappa$ with values in $L|\kappa$ via
such ${\cal X}\to R$ that ${\cal X}\to R$ satisfies the 
properties 0)--6) of subsection A). Among all these possible
choices we consider the ones in which $\dim(R)$ is minimal.

\par

In the notations from the proof of (1) above, and the 
hypothesis of (2), if $\varphi$ is a ${\cal B}$-separable 
approximation, then ${\cal C}=\varphi({\cal B})$ is a 
separable transcendence basis of $L|\lambda$. Now let 
$\varpi:K\to{\cal Y}$ be a ${\cal C}$-separable approximation 
of $L|\lambda$ with values in $K|\kappa$ via ${\cal Y}\to S$. 
(Such approximations do exist, as $L|\lambda$ is supposed 
to have enough separable approximations with values in 
$K|\kappa$.) As in the proof of assertion~(1), there 
exists a projective model ${\cal X}'\to R'$ of $K|\kappa$, 
and a $k$-homomorphism $q:S\to R'$ together with an 
homomorphism of $R'$-varieties $\psi_{R'}:{\cal X}'\to
{\cal Y}_{R'}$, such that $\varpi$ is the base change
of $\psi_{R'}$ to $\kappa$ via the inclusion 
$R'\hookrightarrow\kappa$. Putting everything together,
we obtain a smooth, ${\cal B}$-separable approximation 
of $K|\kappa$ with values in $K|\kappa$ which is defined 
via the composition $p':R\hor{p}S\hor{q}R'$ as follows 
$$
\varphi':{\rm Spec}\,K\>\hookrightarrow\>
     {\cal X}'\hhb5\horr{\phi_{R'}}\hhb5{\cal Y}_{R'}\hhb5
         \horr{{\rm can}}\hhb5{\cal Y}\hhb5\horr{\phi_{S}}
             \hhb5{\cal X}_S\hhb5\horr{{\rm can}}\hhb5{\cal X}
$$
In the above context, let $\kappa^\nix_0={\rm Quot}(R)$, and 
$K_0|\kappa^\nix_0$ be the function field of ${\cal X}\to R$.
Thus $K=K_0\,\kappa$ is the compositum of $K_0$ and $\kappa$ 
over $\kappa^\nix_0$ inside $K$. We now have:

\msn
 
{\it Claim 1.\/} \ $p'$ is injective, and the
$\kappa$-homomorphism $\imath^\nix_{\varphi'}:
\kappa(X_{\varpi'})\to K$ is an isomorphism. Equivalently, 
$K=\imath^\nix_{\varphi'}(K_0)\>\kappa$ inside $K$.
\msn

{\it Proof of Claim 1.\/} \ Let $X'={\cal X}'\times_{R'}\kappa$,  
and $\phi_{\varphi'}:X'\to X_{\varphi'}$ the base change to 
$\kappa$ of the $R'$-morphism ${\cal X}'\to{\cal X}_{R'}$ 
from above. Since both $\varphi$ and $\varpi$ were --among 
other things-- separable approximations, it follows that 
$\phi_{\varphi'}$ is a generically finite and separable 
$\kappa$-morphism; and it defines the above homomorphism of 
function fields $\imath^\nix_{\varphi'}:\kappa(X_{\varpi'})\to K$. 
Let $X\to\kappa$ be the generic fiber of ${\cal X}\to R$. 
Then $X$ is birationally equivalent to $X'$. Therefore,  
$\imath^\nix_{\varphi'}$ is defined by some dominant 
generically finite and separable rational $\kappa$-map 
$\phi':X\birmap X_{\varpi'}$. Thus we have the following 
situation:
\ssn
\item{-} Both $X\to\kappa$ and $X_{\varpi'}\to\kappa$ 
are projective smooth varieties of general type, and 
$\phi':X\birmap X_{\varpi'}$ is a generically finite and 
separable rational $\kappa$-map.
\ssn
\item{-} Let $\omega^\nix_X$ and $\omega^\nix_{X_{\varpi'}}$ be 
their canonical classes. Then for all $m>0$ some given range 
$m\leq n$, we have $\dim|\omega^{\otimes m}_X|=\dim|\omega^
{\otimes m}_{X_{\varpi'}}|=:N_m$.

\ssn
We now choose from the beginning $n$ sufficiently large, such 
as to obatain an $n$-canonical embeddings in ${\lv P}^N_\kappa$,
where $N=N_n$. Thus the same is true for a model ${\cal X}\to R$ 
of $X\to\kappa$ which satisfies the conditions 0)--6) of 
subsection~A) above. In particular, the same is true for 
$X_{\varphi'}$, as it is a base change of fibers of 
${\cal X}\to R$.

\par

Now using e.g.\ Iitaka [I], Ch.5, \S5.4, it follows that in
this situation, every dominant rational separable $\kappa$-map 
$\phi':X\birmap X_{\varpi'}$ is birational.  

\par

In order to finish the proof of the Claim, we remark 
that denoting $\overline R=p'(R)$ the image of $R$ 
under $p'$ (which is the same as the image of $R$ 
under $\imath^\nix_{\varphi'}$), it follows that 
$\overline{\cal X}={\cal X}\times_R\overline R$ has
the property: $\overline{\cal X}\times_{\overline R}\kappa=
X_{\varphi'}$. Since $\kappa(X_{\varphi'})=K$, it follows 
that $\overline{\cal X}\to\overline R$ is a projective
smooth model of $K|\kappa$. Further, $\overline{\cal B}=
\imath^\nix_{\varphi'}({\cal B})$ is defined over $\overline R$,
and $\overline{\cal X}={\cal X}\times_R\overline R$ endowed
with $\overline{\cal B}$ has the properties 0)--6) from 
subsection A) above. By the minimality of $\dim(R)$, it 
follows that $\dim(R)=\dim(\overline R)$, thus $p'$ maps 
$R$ isomorphically onto $\overline R=p'(R)$.

\ssn

Now let ${\cal X}\to R$ a projective smooth model of 
$K|\kappa$ endowed with a separable transcendence basis
${\cal B}$, such that $\dim(R)$ is minimal, and the 
conditions 0)--6) of subsection~A above are satisfied.

\ssn

{\it Claim 2.\/} Let $\varphi:L\to{\cal X}$ be a 
${\cal B}$-separable approximation of $K|\kappa$ via
${\cal X}\to R$. Then $\varphi$ is dominant, and the 
corresponding $\imath^\nix_{\varphi}:\lambda(X_\varphi)\to L$ 
is an isomorphism. In particular, $L|\lambda$ is a function 
field of general type, and $X_\varphi$ is a projective
smooth model of $L|\lambda$.
\ssn
 
{\it Proof of Claim 2.\/} Let us denote ${\cal C}
=\imath^\nix_\varphi({\cal B})$. Then ${\cal C}$ 
is a (separable) transcendece basis of $L|\lambda$. 
Moreover, since ${\cal X}\to R$ has property~2) of
subsection~A), it follows that $[K:\kappa({\cal B})]=
[\lambda(X_\varphi):\lambda({\cal C})]$. Thus
$\imath^\nix_\varphi:\lambda(X_\varphi)\to L$ is an
isomorphism \eqv \ 
$[K:\kappa({\cal B})]=[L:\lambda({\cal C})]$. By Claim~1
(and in the notations from their), it follows that 
$p':R\to S\to R'$ is injective, thus $p:R\to S$ is 
injective. Equivalently, $\varphi:L\to{\cal X}$ 
factoriezes through the generic point of $R$. In other 
words, denoting $\kappa_0={\rm Quot}(R)$, and by $K_0$ 
the function field of ${\cal X}$ we have:
\ssn
\item{-} $K=K_0\,\kappa$ in a canonical way. 
\ssn
\item{-} $\imath^\nix_\varphi:K_0|\kappa_0\to L|\lambda$
is a morphism of function fields.
\ssn
\item{-} ${\cal C}=\imath^\nix_\varphi({\cal B})$ is a 
separable transcendence basis of $L|\lambda$.

\ssn
Now in the notations from the proof of assertion (1) above,
suppose that the model ${\cal Y}\to S$ endowed with the
separable transcendence basis ${\cal C}$ of $L|\lambda$ 
satisfies condition 2) of subsection~A). Then ${\cal B}':=
\imath^\nix_\varpi({\cal C})$ is a separable transcendence 
basis of $K|\kappa$ such that $[\kappa(X_\varpi):\kappa({\cal B}')]
=[L:\lambda({\cal C})]$.
\ssn
Let us set $L_0|\lambda_0=\imath^\nix_\varphi(K_0|\kappa_0)$.
By Claim~1, on has $K=\kappa(X_\varpi)$, and therefore also 
$[K:\kappa({\cal B}')]=[K:\kappa({\cal B})]$. Hence we have 
$[K:\kappa({\cal B})]=[L:\lambda({\cal C})]$, or equivalently, 
$L=L_0\,\lambda$. 

\msn
The proof of assertion~(2) is completed. Thus Theorem~2.1 is proved.
\gata

\bsn

\centerline{\bf 3.\ Proof of Theorem A and Theorem B}

\bsn
Both in the context of Theorem~A and Theorem~B, let 
$\Psi:\ust{K}\to\ust{L}$ be a fixed  isomorphism of some 
ultra-powers $\ust{K}=K^I/{\cal D}$ and $\ust L=L^J/{\cal E}$ 
of $K$ respectively $L$. We set $\ust\kappa$ and $\ust\lambda$ 
for the corresponding ultra-powers of $\kappa$ and $\lambda$ 
inside $\ust K$, respectively $\ust L$. We think about $K|\kappa$ 
and $L|\lambda$ as being diagonally emebdded in $\ust K$, 
respectively $\ust L$. The following holds:

\msn
{\bf Fact 3.1.} \ {\it In the above kontext one has:
\ssn
\item{\rm(1)} If $\kappa$ is a number field or a finite 
field, then $\Psi(\kappa)=\lambda$.
\ssn
\item{\rm(2)} If $\kappa$ is algebraiclly closed, then 
the same is true for $\lambda$, and $\Psi(\ust\kappa)=
\ust\lambda$.
\ssn
\/}

\proof. The first assertion is clear. The second might be 
also known, but we cannot give a reference. It follows 
immediately from Lemma~3.3 below (which itself relies on
the next Lemma~3.2).

\par

We begin by a little preparation. Let ${\lv P}^1_t$ be the 
$t$-projective line over $\kappa$, and $K_0=\kappa(t)$ 
its function field. For every $a\in\kappa$ we denote by 
$K_a\mid K_0$ a minimal field extension of $K_0$ in which 
$\pta{t}\in K_0[T]$ has a root, where
\ssn
\item{a)} $\pta t=T^2-1/(t-a)$, \ if $\chr(\kappa)\neq2$
\ssn
\item{b)} $\pta t=T^2-T-1/(t-a)$, \ if $\chr(\kappa)=2$
\smallskip
Let $C_a\to{\lv P}^1_t$ be the normalisation of ${\lv P}^1_t$
in the Galois field extension $K_a\mid K_0$. It is clear that
$C_a\to{\lv P}^1_t$ is ramified exactly in $t=a,\infty$ in 
case~a), respectively in $t=a$ in case~b). Thus, if $S\subset
\kappa$ is a finite subset of cardinality $n>0$, then the field 
extensions $K_a|K_0$ are linearly disjoint over $K_0$, and their 
compositum $K_S|K_0$ is the function field of the fiber product 
$C_S$ of all the $C_a$ over ${\lv P}^1_t$. Using the Hurwitz 
genus formula, we se that the genus $g_n$ of $C_S$ is given by 
\ssn
\item{a)} $g_n=2^{n-2}(n-3)+1$, \ if $\chr(\kappa)\neq2$ 
\ssn
\item{b)} $g_n=2^{n-1}(n-2)+1$, \ if $\chr(\kappa)=2$,    
\ssn
thus it depends only on $n$ and on $\chr(\kappa)$. The main 
technical point (which on the other hand might be well known
to specialists, but again, we cannot give a reference) is
the following:
\msn
{\bf Lemma 3.2.} \ {\it Let $K|\kappa$ be a function field
with $\kappa$ algebraically closed. For $x\in K$ let $S_x$ 
be the set of all $a\in\kappa$, $a\neq x$, such that $\pta x$
has a root in $K$. Then there exits a bound $c=c^\nix_{K|\kappa}$ 
with the property: If $x\in K$ is a non-constant function, then
the cardinality of $S_x$ is bounded by $c$. To the contrary, if 
$x\in\kappa$, then $S_x=\kappa\backslash\{a\}$, thus infinite.\/}
\msn

{\it Proof of Lemma 3.2.\/} We make induction on $d=\td(K|\kappa)$.
Let $X\to\kappa$ be a projective normal model of $K|\kappa$,
and let $\imath:X\hookrightarrow{\lv P}^N_\kappa$ be a projective 
embedding of $X$. Then for every $d-1$ hyper-planes $H_i$ in 
``general position'' in ${\lv P}^N_\kappa$, let $C=X\cap_i H_i$ 
be the resulting curve in $X$. It is well known that the following
holds:
\ssn
\item{-} $C\to\kappa$ is a normal curve.
\ssn
\item{-} The set of generic points $\eta^\nix_C$ of all 
the generic curves $C$ is dense in $X$.
\ssn
\item{-} The genus $g$ of $C$ is independent of the concrete 
choice of the hyper-planes, it being an invariant of 
the projective embedding $\imath$.

\ssn
We will show that we can take $c=c^\nix_{K|\kappa}=g+2$. 

\par

For every non-constant function $x\in K$, let $K_x$ be the 
relative algebraic closure of $\kappa(x)$ in $K$. Further, 
let $C_x\to\kappa$ be a projective normal (thus smooth) model 
of $K_x|\kappa$. For a be a finite subset $S$ of $S_x$, let 
$K_S\subset K_x$ be the extension of $\kappa(x)$ generated 
by the roots of all the $\pta x$ with $a\in S$. Thus if $C_S$ 
is a projective normal model for $K_S$, then there exists a 
dominant $\kappa$-homomorphism $C_x\to C_S$. Furthermore, 
if $X\to\kappa$ is a normal model for $K|\kappa$, then the
inclusion $K_x\to K$ is defined by a dominant rational
$\kappa$-map $f:X\birmap C_x$.

\ssn

First suppose $\td(K|\kappa)=1$. Then $C_x=X$ is a projective 
normal model for $K|\kappa$, and $g$ is the genus of $X$.
Since $X=C_x\to C_S$ is dominant, by the Hurwitz genus formula, 
$g\geq g_n$, where $n=|S|$ is the cardinality of $S$. 
Since $g_n+2\geq n$, we deduce that all finite subsets 
$S$ of $S_x$ have cardinality bounded by $g+2$.  

\ssn

Coming to the general case, since the set of the generic 
points $\eta^\nix_C$ is dense in $X$, it follows that there 
exit points $\eta^\nix_C$ at which $f:X\birmap C_x$ is 
defined. This means that $f$ defines a rational $\kappa$-map
$f_x:C\birmap C_x$ of projective normal curves. Again, by the 
the Hurwitz genus formula, it follows that $g\geq g_x$ with
$g_x$ the genus of $C_x$. We conclude as above.

\ssn
 
The proof of the Lemma~3.2 is finished.~\gata

\ssn
\msn
{\bf Lemma 3.3.} \ {\it Let $K|\kappa$ be a function 
with $\kappa$ an algebraically closed field. Let further 
$\ust K=K^I/{\cal D}$ be some ultra-power of $K$. Then 
$\ust\kappa$ is the unique maximal algebraically 
closed subfield of $\ust K$.
\smallskip
More precisely, the description of $\ust\kappa$ inside $\ust K$
can be given as follows: For every $\ust x\in\ust K$ let 
$S_{\ust x}\subset\abs\kappa$ be the set of all absolute
algebraic elements of $\kappa$ for which $\pta x$ has roots 
in $\ust K$. Then one has: 
\ssn
\item{\rm(i)} $|S_{\ust x}|\leq c^\nix_{K|\kappa}$ for all 
$\ust x \notin\ust\kappa$, in particular $|S_{\ust x}|$ is finite.
\ssn
\item{\rm(ii)} $|S_{\ust x}|=\abs\kappa\backslash\{\ust x\}$ for all 
$\ust x \in \ust\kappa$, in particular $|S_{\ust x}|$ is infinite.
\ssn\/} 

{\it Proof of Lemma~3.3\/} immediately follows from Lemma~3.2,
thus we will omit the proof here.

\ssn

Finally, the proof of the second assertion of Fact~3.1
follows from Lemma~3.3 and the observation that 
$\Psi(\abs\kappa)=\abs\lambda$ inside $\Psi(\ust K)=\ust L$. 

\bsn
We now come to the {\it Proof of Theorem~A and Theorem~B.\/}

\msn
We first remark that assertion (1) of both Theorem~A
and Theorem~B follows from Theorem~1.4 using Fact~3.1
above. Indeed, suppose first that $\chr(\kappa)=p$ is 
positive. Then $\chr(\lambda)=p$ is positive. Further,
Theorem~1.4, (1), asserts: $\td(K|\kappa)=d$ is equivalent
to saying that $\exists{\cal B}_d=(t_1,\dots,t_d)$ which
is $p$-base of $K$. Moreover, if this assertion is true
for ${\cal B}_d$, then ${\cal B}_d$ is a separable 
transcendence basis of $K|\kappa$. The case $\kappa$ is
algebraically closed of characteristic zero is completely
similar. Now suppose that $\kappa$ is a number field.
Then $\lambda$ is a number field too, as $\kappa=
\abs{\ust\kappa}$, and $\lambda=\abs{\ust\lambda}$,
and $\Psi(\abs{\ust\kappa})=\abs{\ust\lambda}$. Moreover,
using Theorem~1.4, (3), it follows that $\td(K|\kappa)=
\td(L|\lambda)$. 

\msn

We now come to the proof of assertion~(2) of
Theorem~A and Theorem~B. We will be using Theorem~2.1.
In order to do this, we will show that -roughly speaking-
{\it elementary equivalence implies existence of enough
approximations,\/} even without the supplementary 
condition on $K|\kappa$ of being of general type.

\msn
{\bf Key Lemma.} \ {\it Let $K|\kappa$ and $L|\lambda$
be function fields over either finite fields, or number
fields, or algebraically closed fields. Suppose that $K$ 
and $L$ are \ee. Then $K|\kappa$ has enough separable
approximations with values in $L|\lambda$, and $L|\lambda$ 
has enough separable approximations with values in $K|\kappa$.\/}
\msn

{\it Proof.\/} Let $\Psi:\ust K\to \ust L$ be an
isomorphism of some ultra-powers of $K$ and $L$. We
will show that $\Psi$ gives rise to enough approximations
of $K|\kappa$ with values in $L|\lambda$. 

\par

Let ${\cal B}=(t_1,\dots,t_d)$ be a separable transcendence
basis of $K|\kappa$, and further ${\cal X}={\rm Proj}\Big
(R[\underline X]/(\underline f)\Big)\to R$ some model of
$K|\kappa$, where $R=k[\underline y]$ is a finitely generated 
$k$-algebra, $k$ being the prime field of $\kappa$ and 
$\lambda$. We will suppose that ${\cal X}\to R$ satisfies/has 
the properties 0)--4) of Section 2,~A), respectively 
properties 0)--6) if it is the case.  

\par

Without loss of generality, suppose that all the polynomials 
in the system $\underline f$ are absolutely irreducible, and 
$X_0$ is not among them. Then we have $K=\kappa(\underline u)$, 
where $\underline u=(x_1/x_0,\dots,x_n/x_0)$ is the ``canonical 
system'' of rational functions on ${\cal X}$ generating $K|\kappa$. 

\smallskip

Via the the restriction $p=\Psi|_R$, we get an embedding
$p:R\to\ust\lambda$. Unfortunately we do not know {\it a 
priori\/} that $S:=p(R)\subset\lambda$... Via $p$ we can 
consider the base change $\ust{\cal X}={\cal X}\times_R S
\to S\subset\ust\lambda$. This is nothing but $\ust{\cal X}=
{\rm Proj}\Big(S[\underline X]/p(\underline f)\Big)\to S$. 

\par

Next, setting $R=k[\underline y]$, and taking local
representatives $\underline y_j\in\lambda$ for 
$p(\underline y)$, it follows that locally we get local 
representatives $p_j:R\to S_j\subset\lambda$ for the 
ring homomorphism $p:R\to S\subset\ust\lambda$.

\par

Considering ``local representatives'' $\underline f_j$ for 
$p(\underline f)$, we get local representatives ${\cal X}_j
\to S_j$ for $\ust{\cal X}\to S\subset\ust\lambda$, which 
are defined by 
$$
{\cal X}_j=
   {\rm Proj}\Big(S_j[\underline X]/(\underline f_j)\Big)\to S_j.
$$
We remark that ${\cal X}_j\to S_j$ is nothing but the base 
change ${\cal X}\times_RS_j\to S_j$ of ${\cal X}\to R$ to
$S_j$ via the local morphisms $p_j:R\to S_j$. Further, by the 
usual properties of ultra-products we have: If $\underline u_j$
as systems of elements in $L$ are local representatives for 
$\Psi(\underline u)$ as systems of elements of $\ust L$, then 
locally $(1,\underline u_j)$ is a zero of the the homogeneous 
system of equations $\underline f_j=0$ over $S_j$. Moreover, 
for $1\leq\nu\leq d$ let $t_\nu=g_\nu(\underline x)/h_\nu
(\underline x)$ be fixed representations of $t_\nu$ with 
$g_\nu(\underline X)$ and $h_\nu(\underline X)$ homogenous 
equal degree polynomials over $R$. Then $g_{j\nu}=
p_j(g_\nu)$ and $h_{j\nu}=p_j(h_\nu)$ are homogenous equal 
degree polynomials over $S_j$ which locally are representatives
for $p(g_\nu)$, respectively $p(h_\nu)$, and further holds: 
$t_{j\nu}:=g_{j\nu}(\underline x_j)/h_{j\nu}(\underline x_j)$ 
are defined locally as elements of $L$, and they are local 
representatives for $\Psi(t_\nu)$. We denote ${\cal B}_j=
(t_{j1},\dots,t_{jd})$.

\ssn

Now suppose that ${\cal B}_d:={\cal B}$ has supplementary 
properties, as indicated in Theorem~1.4. More precisely,
suppose that:
\ssn
\item{a)} If $\kappa$ and thus $\lambda$ have positive
characteristic $p>0$, then ${\cal B}$ is a separable 
transcendence basis of $K|\kappa$. 
\ssn

\item{$\bullet$} Then the form $q^{(p)}_{\cal B}$ from loc.cit.\ 
does not represent $0$ over $K$, thus over $\ust K$. Therefore, 
$\Psi(q^{(p)}_{\cal B})=q^{(p)}_{\Psi({\cal B})}$ does not 
represent $0$ over $\ust L$. Equivalently, the form 
$q^{(p)}_{{\cal B}_j}$ does not represent $0$ over $L$ locally. 
Again, by Theorem~1.4, (1), it then follows that ${\cal B}_j$ 
is locally a separable transcendence basis of $L|\lambda$.
\ssn
\item{b)} If $\kappa$ is algebraically closed of 
characteristic zero, thus the same is true for $\lambda$, 
then the Pfister form $q^\nix_{\cal B}$ does not 
represent $0$ over $K$, thus over $\ust K$. 
\ssn

\item{$\bullet$} Then $\Psi(q^\nix_{\cal B})=q^\nix_{\Psi({\cal B})}$ 
does not represent $0$ over $\ust L$. Equivalently, the 
Pfister form $q^\nix_{{\cal B}_j}$ does not locally represent 
$0$ over $L$. Again, by Theorem~1.4, (2), it follows that 
${\cal B}_j$ is locally a separable transcendence basis of 
$L|\lambda$.
\ssn
\item{c)} If $\kappa$ be a number field, thus
$\lambda$ is a number field too, then there exit $t_{d+1},
t_{d+2}\in\kappa$ such that the corresponding Pfister form 
$q^\nix_{(t_1,\dots,t_{d+2})}$ does not represent $0$ over 
$K[\sqrt{-1}]$, thus over $\ust K[\sqrt{-1}]$. 
\ssn

\item{$\bullet$} Reasoning as above, and taking into account that 
$u_{d+1}=\Psi(t_{d+1})$ and $u_{d+2}=\Psi(t_{d+2})$ are in 
$\lambda$, it follows that the corresponding Pfister form 
$q_{(t_{j1},\dots,t_{jd},u_{d+1},u_{d+2})}$ does not 
represent $0$ over $L[\sqrt{-1}]$. By Theorem~1.4,~(3), it 
follows that ${\cal B}_j$ is locally a transcendence basis of 
$L|\lambda$. 

\ssn

From the discussion above we deduce that the subfield 
$L_j=\lambda(\underline u_j)$ of $L$ generated by 
$\underline u_j$ locally has transcendence degree $d$ over 
$\lambda$, thus $\underline u_j$ is a generic point of the 
variety ${\cal X}_j\to S_j$. Equivalently, we have obtained 
an approximation 
$$
\varphi_j:{\rm Spec}\,L\to{\cal X}_j\hor{{\rm can}}{\cal X}
$$
where ${\cal X}_j\hor{{\rm can}}{\cal X}$ comes from the base 
change ${\cal X}_j={\cal X}\times_R S_j$. By its construction, 
the approximation $\varphi_j$ is separable, provided ${\cal B}$ 
a separable transcendence base.

\msn

One concludes the proof of Theorem~A and Theorem~B
by applying Theorem~2.1.

\msn
{\it Remark.\/} \ Under the hypothesis of assertion~(2) 
of Theorem~A and Theorem~B, it follows that $\Psi(K\ust\kappa)=
L\ust\lambda$.  

\bsn
 
\AuthorRefNames[F--H--V]

\references

[B--T]
\name{Bass H.} and \name{Tate, J.}, The Milnor ring of a global field. 
    Algebraic $K$-theory, II: "Classical" algebraic $K$-theory and 
    connections with arithmetic (Proc.\ Conf., Seattle, Wash., Battelle 
    Memorial Inst., 1972), pp. 349--446. LNM {\bf 342}, Springer, Berlin 
    1973.

[B--S]
\name{Bell J.-L.} and \name{Slomson, A.-B}, Models and ultra-products:
an introduction, Amsterdam 1969.

[B]
\name{Bloch, S.}, .

[BOU]
\name{Bourbaki}, {\it Alg\`ebre commutative,\/} Hermann 
     Paris 1964.

[C--S]
\name{Cornell--Silverman}, Introduction to arithmetic geometry,
eds, .

[D\hhb{1}1]
\name{Duret, J.-L.}, Sur la th\'eorie \'el\'ementaire des corps 
de fonctions, J.~Symbolic Logic {\bf 51} (1986), 948--956.

[D\hhb{1}2]
\bibline, \'Equivalence \'el\'ementaire et isomorphisme des corps
de courbe sur un corps algebriquement clos, J.~Symbolic Logic 
{\bf 57} (1992), 808--923.

[D--S]
\name{van den Dries, L.} and \name{Schmidt, K.}, Bounds in the theory
of polynomial rings over fields. A non-standard approach, Invent.~Math.
{\bf 76} (1984), 77--91.

[F]
\name{Faltings, G.}, Endlichkeitss\"atze f\"ur abelsche
     Variet\"aten \"uber Zahlk\"orpern, Invent.\ Math. {\bf 73}
     (1983), 349--366.

[Ei]
\name{Eisenbud, D.}, Commutative Algebra, GTM 150, Springer, 1999.

[E--L]
\name{Elman, R.} and \name{Lam T.-Y.}, Pfister forms and K-theory 
    of fields, J. Algebra {\bf 23} (1972), 181--213. 

[H]
\name{Hartshorne, R.} Algebraic Geometry, Springer Verlag, 
     New York $\cdot$ Berlin $\cdot$ Heidel\-berg . . . 1993 
     (sixth printing).

[I]
\name{Iitaka, Sh.} Algebraic Geometry, Springer Verlag, 
     New York $\cdot$ Berlin $\cdot$ Heidel\-berg 1982.

[J--R]
\name{Jacobs, B.} and \name{Rost, M.}, Degree four cohomological 
     invariants for quadratic forms, Invent.\ Math.\ {\bf } (1989), 
     552--570.

[Kh]
\name{Kahn, E.}, La conjecture de Milnor (d'apr\`es Voevodsky),
      S\'eminaire Bourbaki, Asterisque {\bf 245} (1997), 379--418.

[Ka]
\name{Kani, E.}, Bounds on the number of non-rational subfields of a
          function field, Invent.\ Math.\ {\bf 85} (1986), 185--198.

[Ko]
\name{Koch, H.}, Die Galoissche Theorie der $p$-Erweiterungen,
Math.\ Monographien {\bf 10}, Berlin 1970.

[M--S\hhb{1}1]
\name{Merkurjev, A.-A.} and {Suslin, A.},
Merkurjev, A.~S., ${\rm K}_2$ of a field, and the Brauer group,
          Sov.Math.\ Doklady {\bf 24} (1981), 546--551.

[M--S\hhb{1}2]
\name{Merkurjev, A.~S.} and \name{Suslin A.~A.}, K-cohomology of 
       Brauer--Severi varieties and the norm residue homomorphism, 
       Math.\ USSR Isz.\ {\bf ???} (1983), 307--340.

[M\hhb{1}1]
\name{Milnor, J.}, Algebraic K-theory and quadratic forms, Invent.\ 
      Math.\ {\bf 9} (1970), 318--344.

[M\hhb{1}2]
\bibline, Introduction to algebraic K-theory, Ann.\ Math.\ Studies 
          {\bf 72}, Princeton Univ.\ Press 1971.

[Pa]
\name{Parshin, A. N.}, Finiteness Theorems and Hyperbolic Manifolds,
     in {\it The Groth\-endieck Festschrift\/} III, ed P.~Cartier 
     et all, PM Series Vol 88, Birkh\"auser Boston Basel Berlin 1990.

[Pi]
\name{ Pierce, D.}, Function fields and elementary equivalence,
Bull.\ London Math.\ Soc.{\bf 31} (1999), 431--440.

[P]
\name{ Pop, F.}, On Grothendieck's conjecture of anabelian geometry,
Ann of Math. {\bf 138} (1994), 145--182. 

[R--R]
\name{Robinson, A.} and \name{Roquette, P.}, On the finiteness of
Siegel and Mahler concerning Diophantine equations, J.~Number Theory
{\bf 7} (1975), 121--176. 

[Q]
\name{Quillen, D.}, Higher algebraic K-Theory I, in: Springer LNM 341, 
          85--147. 

[R]
\name{Rost, M.}, see References in [Kh] above.

[S]
\name{Serre, J.-P.}, {\it Cohomologie Galoisienne,\/} LNM 5,
     Springer 1965.

[S--V]
\name{Suslin, A.\ A.} and \name{Vojevodsky, V.}, Bloch-Kato conjecture 
      and motivic cohomology with finite coefficients, preprint (1995).

[V]
\name{Vojevodsky, V.}, see References in [Kh] above.

[Z--S]
\name{Zariski, O.} and \name{Samuel, P.}, {\it Commutative Algebra,\/}
     Van Nostrand, Princeton 1958 and 1960.

\endreferences

\bye